\documentclass[11pt]{article}
\usepackage{epic,latexsym,amssymb}
\usepackage{amsfonts}
\usepackage{amscd}
\usepackage{amsmath}
\usepackage{graphicx}
\usepackage{color}
\usepackage{caption,
caption}
\usepackage{tikz}

\usepackage{float}

\newtheorem{thm}{Theorem}

\newtheorem{ob}[thm]{Observation}
\newtheorem{prop}[thm]{Proposition}
\newtheorem{lem}[thm]{Lemma}

\newtheorem{claim}{Claim}

\newcommand{\trim}{\mathrm{trim}}

\newcommand{\diam}{{\rm diam}}

\newcommand{\pc}{{\rm pc}}

\newcommand{\cP}{\mathcal{P}}
\newcommand{\cQ}{\mathcal{Q}}

\newcommand{\cT}{\mathcal{T}}

\newcommand{\proof}{\noindent\textbf{Proof. }}

\newcommand{\smallqed}{{\tiny ($\Box$)}}
\newcommand{\qed}{$\Box$}

\newcommand{\1}{\vspace{0.1cm}}

\newenvironment{unnumbered}[1]{\trivlist
\item [\hskip \labelsep {\bf #1}]\ignorespaces\it}{\endtrivlist}

\newcommand{\QEDmark}{\mbox{\textsc{qed}}}
\newcommand{\proofStarter}[1]{\textsc{#1} }

\def\vertex(#1){\put(#1){\circle*{2}}}
\def\vertexo(#1){\put(#1){\circle{2}}}
\def\vert(#1){\put(#1){\circle*{1.5}}}
\def\verto(#1){\put(#1){\circle{1.5}}}
\def\lab(#1)#2{\put(#1){\makebox(0,0)[c]{#2}}}
\definecolor{DarkGreen}{rgb}{0.2, 0.6, 0.3}

\definecolor{electricindigo}{rgb}{0.44, 0.0, 1.0}

\let\oldenumerate\enumerate
\renewcommand{\enumerate}{
  \oldenumerate
  \setlength{\itemsep}{0.5pt}
  \setlength{\parskip}{0pt}
  \setlength{\parsep}{0pt}
}

\begin{document}

\title{Matching, Path Covers, and Total Forcing Sets}
\author{$^{1,2}$Randy Davila and $^{1}$Michael A. Henning\\
\\
$^1$Department of Pure and Applied Mathematics\\
University of Johannesburg \\
Auckland Park 2006, South Africa \\
\small {\tt Email: mahenning@uj.ac.za} \\
\\
$^2$Department of Mathematics and Statistics\\
University of Houston--Downtown \\
Houston, TX 77002, USA \\
\small {\tt Email: davilar@uhd.edu}
}

\date{}
\maketitle

\begin{abstract}
A dynamic coloring of the vertices of a graph $G$ starts with an initial subset $S$ of colored vertices, with all remaining vertices being non-colored. At each discrete time interval, a colored vertex with exactly one non-colored neighbor forces this non-colored neighbor to be colored. The initial set $S$ is called a forcing set of $G$ if, by iteratively applying the forcing process, every vertex in $G$ becomes colored. If the initial set $S$ has the added property that it induces a subgraph of $G$ without isolated vertices, then $S$ is called a total forcing set in $G$. The minimum cardinality of a total forcing set in $G$ is its total forcing number, denoted $F_t(G)$. The path cover number of $G$, denoted $\pc(G)$, is the minimum number of vertex disjoint paths such that every vertex belongs to a path in the cover, while the matching number of $G$, denoted $\alpha'(T)$, is the number of edges in a maximum matching of $G$. Let $T$ be a tree of order at least two. We observe that $\pc(T) + 1 \le F_t(T) \le 2\pc(T)$, and we prove that $F_t(T) \le \alpha'(T) + \pc(T)$. Further, we characterize the extremal trees achieving equality in these bounds.  
\end{abstract}


\newpage
\section{Introduction}

Coloring the vertices of a graph $G$ and allowing this initial coloring to propagate throughout the vertex set of $G$ is known as a \emph{dynamic coloring} of $G$. In this paper, we focus on the dynamic coloring due to the \emph{forcing process}, which is defined in~\cite{DaHe17a} as follows: Let $G$ be a finite and simple graph with vertex set $V(G)$, and let $S\subseteq V(G)$ be a set of initially ``colored'' vertices, all remaining vertices being ``uncolored''. All vertices contained in $S$ are said to be $S$-colored, while all vertices not in $S$ are $S$-uncolored. At each discrete time step, if a colored vertex has exactly one uncolored neighbor, then this colored vertex \emph{forces} its uncolored neighbor to become colored. If $v$ is such a colored vertex, then we call $v$ a \emph{forcing vertex}, and say that $v$ has been \emph{played}. The initial set of vertices $S$ is a \emph{zero forcing set}, if by iteratively applying this forcing process all of $V(G)$ becomes colored. Such a set $S$ is called an $S$-\emph{forcing set}. If $S$ is a zero forcing set of $G$ and $v$ is a $S$-colored vertex which has been played, then $v$ is called a \emph{$S$-forcing vertex}. The \emph{zero forcing number} of $G$, written $Z(G)$, is the cardinality of a minimum forcing set in $G$.  The concept of zero forcing in graphs was originally introduced in~\cite{AIM-Workshop} and further studied, for example, in~\cite{AIM-Workshop, k-Forcing, Barioli13, DaHe17a, DaKe15, Genter1, GeRa18, Edholm, Hogben10, Hogben16, LuTang}.

If $S$ is a zero forcing set of $G$ with the additional property that the subgraph of $G$ induced by $S$ contains no isolated vertex, then $S$ is a \emph{total forcing set}, abbreviated TF-set, of $G$. The \emph{total forcing number} of $G$, written $F_t(G)$, is the cardinality of a minimum TF-set in $G$. The notion of total forcing in graphs was first introduced in~\cite{Davila} as a strengthening of zero forcing in graphs and studied further, for example, in~\cite{DaHe17b, DaHe17c, DaHe17d, DaMaSt16}. In this paper, we obtain bounds relating the total forcing number of a tree to its path cover number. Further, we obtain a relationship between the total forcing number, the matching number, and path cover number in a tree.

\medskip
\noindent\textbf{Definitions and Notation.}
For notation and graph terminology, we will typically follow the monograph~\cite{MHAYbookTD}. Specifically, this paper will only consider finite and simple graphs. Let $G$ be a graph with vertex set $V(G)$ and edge set $E(G)$. The \emph{order} and \emph{size} of $G$ will be denoted by $n(G)= |V(G)|$ and $m(G) = |E(G)|$, respectively. Two vertices $u$ and $v$ are \emph{neighbors} in $G$ if they are adjacent, that is, if $uv \in E(G)$. The \emph{open neighborhood} of a vertex $v\in V(G)$, written $N_G(v)$, is the set of all neighbors of $v$, whereas the \emph{closed neighborhood} of $v$ is $N_G[v] = N_G(v)\cup \{v\}$. The \emph{degree} of a vertex $v$ in $G$, written $d_G(v)$, is the number of neighbors of $v$ in $G$; and so, $d_G(v) = |N_G(v)|$. A \emph{nontrivial graph} is a graph of order at least~$2$. We denote the \emph{complete graph}, \emph{path}, and \emph{cycle}, on $n$ vertices by $K_n$, $P_n$, and $C_n$, respectively.

The distance between two vertices $v$ and $w$ in a connected $G$ is the length of a shortest $(v,w)$-path in $G$, and is denoted by $d_G(v,w)$. The maximum distance among all pairs of vertices in $G$ is the \emph{diameter} of $G$, denoted by $\diam(G)$.
A \emph{leaf} is a vertex of degree~$1$, while its neighbor is a \emph{support vertex}. A \emph{strong support vertex} is a vertex with at least two leaf neighbors. A \emph{star} is a non-trivial tree with at most one vertex that is not a leaf. Thus, a star is the tree $K_{1,k}$ for some $k \ge 1$. For $r, s \ge 1$, a \emph{double star} $S(r,s)$ is the tree with exactly two vertices that are not leaves, one of which has $r$ leaf neighbors and the other $s$ leaf neighbors. A \emph{pendant edge} of a graph is an edge incident with a vertex of degree~$1$.

A \emph{rooted tree} $T$ distinguishes one vertex $r$ called the \emph{root}. For each vertex $v \ne r$ of $T$, the \emph{parent} of $v$ is the
neighbor of $v$ on the unique $(r,v)$-path, while a \emph{child} of $v$ is any other neighbor of $v$. The set of children of $v$ is denoted by  $C(v)$. A \emph{descendant} of $v$ is a vertex $u \ne v$ such that the unique $(r,u)$-path contains $v$, while an \emph{ancestor} of $v$ is a vertex $u \ne v$ that belongs to the $(r,v)$-path in $T$. In particular, every child of $v$ is a descendant of~$v$ while the parent of $v$ is an ancestor of $v$. The \emph{grandparent} of $v$ is the ancestor of $v$ at distance~$2$ from $v$. A \emph{grandchild} of $v$ is the descendant of $v$ at distance~$2$ from $v$.
We let $D(v)$ denote the set of descendants of $v$, and we define $D[v] = D(v) \cup \{v\}$.

The \emph{contraction} of an edge $e = xy$ in a graph $G$ is the graph obtained from $G$ by replacing the vertices $x$ and $y$ by a new vertex and joining this new vertex to all vertices that were adjacent to $x$ or $y$ in $G$. Given a non-trivial tree $T$, the \emph{trimmed tree} of $T$, denoted $\trim(T)$, is the tree obtained from $T$ by iteratively contracting edges with one of its incident vertices of degree exactly~$2$ and with the other incident vertex of degree at most~$2$ until no such edge remains. We note that if the original tree $T$ is a path, then $\trim(T)$ is a path $P_2$, while if $T$ is not a path, then every edge in $\trim(T)$ is incident with a vertex of degree at least~$3$. In particular, if $T$ is not a path, then every support vertex in $\trim(T)$ has degree at least~$3$.

A \emph{path cover} of $G$ is a collection of vertex disjoint paths such that every vertex belongs to exactly one path of $G$, and the cardinality of a minimum path cover is known as the \emph{path cover number of} $G$, denoted $\pc(G)$. Path covers are a fundamental concept in graph theory. Papers relating domination parameters and the path cover number can be found, for example, in~\cite{DeLa07,HeWa17}.

Two edges in a graph $G$ are \emph{independent} if they are not adjacent in $G$. A set of pairwise independent edges of $G$ is called a \emph{matching} in $G$, while a matching of maximum cardinality is a \emph{maximum matching}. The number of edges in a maximum matching of $G$ is the \emph{matching number} of $G$ which we denote by $\alpha'(G)$. Matchings in graphs are extensively studied in the literature (see, for example, the classical book on matchings by Lov\'{a}sz and Plummer~\cite{LoPl86}, and the excellent survey articles by Plummer~\cite{Pl03} and Pulleyblank~\cite{Pu95}).

In this paper, we relate the total forcing number of a tree with its path cover number and matching number. We proceed as follows. In Section~\ref{S:main}, we present the statement of our main results. Thereafter, we state some known results in Section~\ref{S:known} that will be helpful in proving our main results. In Section~\ref{S:main1} and Section~\ref{S:main2} we present a proof of Theorem~\ref{t:main1} and Theorem~\ref{t:main2}, respectively. We conclude our discussion in Section~\ref{S:remark} with some remarks  and open problems for future research.

\section{Main Results}
\label{S:main}

We have two immediate aims in this paper. Our first aim is to establish a relationship between the total forcing number of a tree and its path cover number. We shall prove the following result, a proof of which is given in Section~\ref{S:main1}.

\begin{thm}
\label{t:main1}
If $T$ is a nontrivial tree, then
\[
\pc(T) + 1 \le F_t(T) \le 2\pc(T),
\] and all possible values of $F_t(T)$ in this range are possible. Further, the following hold. \\[-24pt]
\begin{enumerate}
\item $F_t(T) = \pc(T) + 1$ if and only if $\trim(T) = P_2$ or $\trim(T) \cong K_{1,n-1}$ for some $n \ge 4$.
\item $F_t(T) = 2\pc(T)$ if and only if $T$ has a unique minimum path cover and every path in this cover starts and ends at distinct leaves of~$T$.
\end{enumerate}
\end{thm}

Our second aim is to establish a relationship between the total forcing number, the matching number and the path cover number of a graph. For this purpose, we define a family of trees $\cT$ as follows. Let $T'$ be an arbitrary tree (possibly, trivial) and let $A$ be a subset of vertices in $T'$ such that either $A = V(T')$ or $V(T') \setminus A$ is an independent set in $T'$ containing no leaf of $T'$. Let $T$ be the tree obtained from $T'$ by attaching at least two pendant edges to each vertex of $A$. We call the tree $T'$ the \emph{underlying tree} of the tree $T$, and we call the set $A$ the \emph{attacher set} of $T$. Further, we call each vertex of $A$ an \emph{attacher vertex} of $T$. We note that the attacher vertices of $T$ are precisely the support vertices of $T$, and each attacher vertex is a strong support vertex of $T$ with all its leaf neighbors outside $T'$. Let $\cT$ be the family of all such trees $T$, together with the tree $K_2$. We shall prove the following result, a proof of which is given in Section~\ref{S:main2}.

\begin{thm}
\label{t:main2}
If $T$ is a nontrivial tree, then  \[F_t(T) \le \alpha'(T) + \pc(T),\] with equality if and only if $T \in \cT$.
\end{thm}

\section{Known Results}
\label{S:known}

The zero forcing number and total forcing number of paths, cycles, complete graphs and stars is easy to compute.

\begin{ob}{\rm (\cite{DaHe17b})}
\label{simpleFormula}
The following holds. \\
\indent {\rm (a)} For $n \ge 2$, $Z(P_n) = 1$ and $F_t(P_n) = 2$. \\
\indent {\rm (b)} For $n \ge 3$, $Z(C_n) = F_t(C_n) = 2$. \\
\indent {\rm (c)} For $n \ge 3$, $Z(K_n) = F_t(K_n) = n - 1$.\\
\indent {\rm (d)} For $n \ge 3$, $Z(K_{1,n-1}) = n-2$ and $F_t(K_{1,n-1}) = n - 1$.
\end{ob}

We recall a useful lemma in~\cite{DaHe17b}.

\begin{lem}{\rm (\cite{DaHe17b})}
\label{leafLem}
If $G$ is an isolate-free graph, then every vertex of $G$ with at least two leaf neighbors is contained in every TF-set, and all except possibly one leaf neighbor of such a vertex is contained in every TF-set.
\end{lem}

The following observation shows that the total forcing number of an isolate-free graph is bounded above by twice the forcing number.

\begin{ob}{\rm (\cite{DaHe17b})}
\label{o:FtZ}
If $G$ is an isolate-free graph, then $F_t(G) \le 2Z(G)$.
\end{ob}

The following results were obtained in~\cite{DaHe17c}.

\begin{lem}{\rm (\cite{DaHe17c})}
\label{lemPaths}
Let $G$ be an isolate-free graph that contains an edge $e$ incident with a vertex of degree at most~$2$. If $G'$ is obtained from $G$ by subdividing the edge $e$ any number of times, then $F_t(G) = F_t(G')$.
\end{lem}

\begin{lem}{\rm (\cite{DaHe17c})}
\label{l:trimT}
If $T$ is a non-trivial tree, then the following hold. \\
\indent {\rm (a)} $F(T) = F(\trim(T))$. \\
\indent {\rm (b)} $F_t(T) = F_t(\trim(T))$. \\
\indent {\rm (c)} The trees $T$ and $\trim(T)$ have the same number of leaves.
\end{lem}

As observed earlier, if $T$ is a non-trivial path, then $\trim(T) = P_2$. The following result establishes a lower bound on the total forcing number in terms of its zero forcing number.

\begin{thm}{\rm (\cite{DaHe17c})}
\label{t:lowerbd}
If $T$ is a non-trivial tree, then $F_t(T) \ge Z(T)+1$, with equality if and only if $\trim(T) = P_2$ or $\trim(T) \cong K_{1,n-1}$ for some $n \ge 4$.
%
%
\end{thm}

The following relation between the zero forcing and path cover numbers was obtained in~\cite{AIM-Workshop,Hogben10}. In particular, the zero forcing number of a tree is precisely its path cover number.

\begin{thm}
\label{t:cover}
The following hold. \\
\indent {\rm (a)}  {\rm (\cite{Hogben10})} If $G$ is a graph, then $Z(G) \ge \pc(G)$. \\
\indent {\rm (b)}  {\rm (\cite{AIM-Workshop,Hogben10})} If $T$ is a tree, then $Z(T) = \pc(T)$.
\end{thm}

\section{Proof of Theorem~\ref{t:main1}}
\label{S:main1}

In this section, we prove Theorem~\ref{t:main1}. For this purpose, we first present a series of preliminary lemmas which will be used in our subsequent argument to establish the desired characterization stated in Theorem~\ref{t:main1}.

\begin{lem}
\label{l:subdivide}
If a tree $T$ contains an edge $e$ with one of its incident vertices of degree exactly~$2$ and with the other incident vertex of degree at most~$2$, then $\pc(T) = \pc(T')$ where $T'$ is obtained from $T$ by contracting the edge $e$. Further, there is a one-to-one correspondence between the minimum path covers in $T$ and $T'$.
\end{lem}
\proof Let $e = uv$, where $d_T(u) \le 2$ and $d_T(v) = 2$. Let $w$ be the neighbor of $v$ different from $u$, and if $d_T(u) = 2$, then let $t$ be the neighbor of $u$ different from $v$. Let $T'$ be obtained from $T$ by contracting the edge $e$, and let $x$ be the resulting new vertex. Thus, in $T'$ either $x$ is a leaf with $w$ as its neighbor or $x$ has degree~$2$ with $t$ and $w$ as its neighbors. Let $\cP$ and $\cP'$ be minimum path covers in $T$ and $T'$, respectively. By the minimality of the path cover $\cP$, the vertices $u$ and $v$ belong to the same path in $\cP$. Let $P_v$ be the path in $\cP$ that contains~$u$ and $v$. Replacing the vertices $u$ and $v$ on $P_v$ with the vertex $x$, and leaving all other paths in $\cP$ unchanged produces a path cover in $T'$, implying that $\pc(T') \le |\cP| = \pc(T)$. Conversely, if $P_x'$ is the path in $\cP'$ that contains the vertex~$x$, then replacing the vertex $x$ on the path $P_x'$ with the deleted vertices $u$ and $v$, and leaving all other paths in $\cP'$ unchanged, produces a path cover in $T$, implying that $\pc(T) \le |\cP'| = \pc(T')$. Thus, $\pc(T) = \pc(T')$ and there is a one-to-one correspondence between the minimum path covers in $T$ and $T'$.~\qed

\medskip
Since every non-trivial tree $T$ can be reconstructed from its trimmed tree $\trim(T)$ by applying a sequence of subdivisions of edges incident with a vertex of degree at most~$2$, as an immediate consequence of Lemma~\ref{l:subdivide} the path cover number of a tree and its trimmed tree are identical. Further, there is a one-to-one correspondence between the minimum path covers in a tree and its trimmed tree. We state this formally as follows.

\begin{lem}
\label{l:trimTb}
If $T$ is a non-trivial tree, then $\pc(T) = \pc(\trim(T))$. Further, there is a one-to-one correspondence between the minimum path covers in $T$ and $\trim(T)$.
\end{lem}

We shall also need the following property of a minimum path cover that contains a strong support vertex.

\begin{lem}
\label{l:strong}
If $v$ is a strong support vertex in a graph $G$ with leaf neighbors $u$ and $w$, then there exists a minimum path cover in $G$ that contains the path $uvw$.
\end{lem}
\proof Let $\cP$ be a minimum path cover in $G$. Let $P_u$, $P_v$ and $P_w$ be the paths in $\cP$ that contain the vertices $u$, $v$ and $w$, respectively. If $P_u = P_v = P_w$, then $P_v$ is the path $uvw$ and we are done. Hence, we may assume renaming $u$ and $w$ if necessary, that $P_u \ne P_v$. Thus, $P_u$ is the trivial path consisting of the vertex~$u$. By the minimality of the path cover $\cP$, the vertex $v$ is an internal vertex of $P_v$.
If $P_v = P_w$, then we replace the two paths $P_u$ and $P_v$ in $\cP$ with the following two paths: the path $uvw$ and the path obtained from $P_v$ by deleting from it the vertices $v$ and $w$, and we leave all other paths in $\cP$ unchanged. If $P_v \ne P_w$, then the path $P_w$ is the trivial path consisting of the vertex~$w$. In this case, we replace the three paths $P_u$, $P_v$ and $P_w$ in $\cP$ with the following three paths: the path $uvw$ and the two paths obtained from $P_v$ by deleting the vertex~$v$, and we leave all other paths in $\cP$ unchanged. Let $\cP'$ denote the resulting new path cover. In both cases, $|\cP'| = |\cP|$, and so $\cP'$ is a minimum path cover in $G$ that contains the path $uvw$.~\qed

\begin{lem}
\label{l:good_cover}
If a tree $T$ has a unique minimum path cover and every path in this cover starts and ends at distinct leaves of~$T$, then the following hold.
\\[-22pt]
\begin{enumerate}
\item The set consisting of a leaf and it neighbor from each path in the minimum path cover forms a TF-set in $T$.
\item $F_t(T) = 2\pc(T)$.
\end{enumerate}
\end{lem}
\proof Let $T$ be a tree and suppose that $T$ has a unique minimum path cover and every path in this cover starts and ends at distinct leaves of~$T$. We proceed by induction on $\pc(T)$ to show that $F_t(T) = 2\pc(T)$. If $\pc(T) = 1$, then $T$ is a path on at least two vertices, and so $F_t(T) = 2 = 2\pc(T)$. Further, the set consisting of a leaf of $T$ and it neighbor forms a TF-set in $T$. This establishes the base case. Let $k \ge 2$ and assume that if $T'$ is a tree with $\pc(T') < k$ that has a unique minimum path cover $\cP'$ and every path in $\cP'$ starts and ends at distinct leaves of~$T'$, then $F_t(T') = 2\pc(T')$ and the set consisting of a leaf and it neighbor from each path in $\cP'$ forms a TF-set in $T$.

Let $T$ be a tree with $\pc(T) = k$ and suppose that $T$ has a unique minimum path cover $\cP$ and every path in this cover starts and ends at distinct leaves of~$T$. Let $T_\cP$ be a graph of order~$k$ whose vertices correspond to the $k$ paths in $\cP$ and where two vertices in $T_\cP$ are joined by an edge if and only if there is an edge between the corresponding paths in $\cP$. Since $T$ is a tree, the graph $T_\cP$ is a tree. Let $v'$ be a leaf in $T_\cP$ and let $v$ be its neighbor in $T_\cP$, and let $P'$ and $P$ be the paths in $\cP$ corresponding to the vertices $v'$ and $v$. We note that both $P'$ and $P$ start and end at distinct leaves of $T$. Since $vv'$ is an edge of $T_\cP$, there is an edge $e = xx'$ that joins an internal vertex $x$ of $P$ and an internal vertex $x'$ of $P'$. Let $\cP' = \cP \setminus \{P'\}$ and let $T'$ be the tree obtained from $T$ by deleting the vertices on the path $P'$; that is, $T' = T - V(P')$.

Since $\cP'$ is a path cover of $T'$, we note that $\pc(T') \le |\cP'| = |\cP| - 1 = \pc(T) - 1$. Every path cover in $T'$ can be extended to a path cover in $T$ by adding to it the path $P'$, implying that $\pc(T) \le \pc(T') + 1$. Consequently, $\pc(T') = \pc(T) - 1 < k$. Thus, $\cP'$ is a minimum path cover in $T'$. If $T'$ has a minimum path cover different from $\cP'$, then such a path cover can be extended to a minimum path cover in $T$ by adding to it the path $P'$ to produce a minimum path cover different from $\cP$, contradicting the fact that $\cP$ is the unique minimum path cover of $T$. Hence, $T'$ has a unique minimum path cover, namely $\cP'$. Every leaf of $T$ that does not belong to the path $P'$ is a leaf of $T'$, and every leaf of $T'$ is a leaf of $T$. Every path in $\cP'$ therefore starts and ends at distinct leaves of~$T'$. Thus, $T'$ is a tree with $\pc(T') < k$ that has a unique minimum path cover and every path in this cover starts and ends at distinct leaves of~$T'$. Applying the inductive hypothesis to $T'$, $F_t(T') = 2\pc(T')$ and the set, $S'$ say, consisting of a leaf and it neighbor from each path in $\cP'$ forms a TF-set in $T'$.

Let $P'$ be a $(u',v')$-path given by $u_1u_2 \ldots u_\ell$ where $u' = u_1$ and $v' = u_\ell$, and so the path $P'$ starts at the leaf $u'$ and ends at the leaf $v'$. Recall that exactly one vertex of $P'$, namely the vertex $x'$, is adjacent in $T$ to a vertex outside $P'$, namely to the vertex $x$ which belongs to the path $P$. Further, $x$ and $x'$ are internal vertices of $P$ and $P'$, respectively. Let $S = S' \cup \{u_1,u_2\}$. We show that $S$ is a TF-set of $T$. Let $x' = u_j$, where  we note that $j \in [\ell - 1] \setminus \{1\}$. If $j > 2$, then $u_2 \ne x'$ and as the first vertex played in the forcing process we play the vertex $u_2$, thereby coloring $u_3$. Further if $j > 3$, then $u_3 \ne x'$ and as the second vertex played in the forcing process we play the vertex $u_3$, thereby coloring $u_4$. Continuing in this way, we play as the first few vertices in the forcing process the vertices $u_2, \ldots, u_{j-1}$, thereby coloring all vertices $u_1,u_2, \ldots, u_j$. Thereafter we play the identical sequence of vertices in the forcing process in $T'$ starting with the set $S'$ that results in all $V(T')$ colored. Since $S'$ is a TF-set of $T'$ and since $x' = u_j$ is colored, we note that this results in all vertices of $T'$ colored. Finally, we play the sequence of vertices $u_j,\ldots,u_{\ell - 1}$ in turn, resulting in all vertices of $P'$ colored. Thus, $S$ is a TF-set of $T$ that consists of a leaf and it neighbor from each path in $\cP$. This proves Part~(a).

We show next that $F_t(T) = F_t(T') + 2$. Let $T^*$ be the tree obtained from $T'$ by adding the path $u'x'v'$ and the edge $xx'$. If $T \ne T^*$, then the path $P'$ has order at least~$4$ and $T'$ can be obtained from $T^*$ by a sequence of edge subdivisions where each edge that is subdivided is incident with a vertex of degree at most~$2$. In this case, Lemma~\ref{lemPaths} implies that $F_t(T) = F_t(T^*)$. If $T = T^*$, then trivially $F_t(T) = F_t(T^*)$. Hence, it suffices for us to show that $F_t(T^*) = F_t(T') + 2$. Let $S$ be a minimum TF-set in $T^*$. By Lemma~\ref{leafLem}, the set $S$ contains the vertex $x'$ and at least one of $u'$ and $v'$. If both $u'$ and $v'$ belong to $S$, then $(S \setminus \{u'\}) \cup \{x\}$ is a minimum TF-set of $T^*$. Hence, we may choose $S$ so that $u' \notin S$ and $\{x',v'\} \subset S$. Thus since $S$ is a TF-set of $T^*$, the set $S \setminus \{x',v'\}$ is a TF-set of $T'$, and so $F_t(T') \le |S| - 2 = F_t(T^*) - 2$. Conversely, every minimum TF-set of $T'$ can be extended to a TF-set of $T^*$ by adding to it the vertices $x'$ and $v'$, and so $F_t(T^*) \le F_t(T') + 2$. Consequently, $F_t(T^*) = F_t(T) - 2$.
By our earlier observations, $F_t(T') = 2\pc(T')$, $\pc(T') = \pc(T) - 1$ and $F_t(T) = F_t(T') + 2$, implying that $F_t(T) = 2\pc(T)$. This proves Part~(b) and completes the proof of Lemma~\ref{l:good_cover}.~\qed

\medskip
We show next that all possible values of $F_t(T)$ in the range from $\pc(T) + 1$ to $2\pc(T)$ are possible.

\begin{prop}
\label{p:range}
For any two given positive integers $k$ and $\ell$ where $k \in [\ell]$, there exists a tree $T$ satisfying $\pc(T) = \ell$ and $F_t(T) = k + \ell$.
\end{prop}
\proof Let $\ell \ge 1$ be an arbitrary integer. If $k = 1$, then taking $T = K_{1,\ell + 1}$ we note that $\pc(T) = \ell$ and $F_t(T) = \ell + 1$, and so $F_t(T) = k + \ell$. If $k = \ell$, then let $T$ be obtained from a path $P_{\ell}$ on $\ell$ vertices by adding two pendant edges to each vertex of the path. The resulting tree $T$ satisfies $\pc(T) = \ell$ and $F_t(T) = 2\ell$, and so $F_t(T) = k + \ell$. Hence, we may assume that $\ell \ge 3$ and $k \in [\ell - 1] \setminus \{1\}$. For integers $r \ge 3$ and $s \ge 1$, let $T' \cong P_{s+1}$ be a path on $s+1$ vertices, and let $T$ be the graph obtained from $T'$ by adding $r$ pendant edges to one vertex of $T'$ and adding two pendant edges to each of the remaining $s$ vertices of $T'$. Thus, $T$ has order~$r + 3s + 1$ and every vertex in $V(T')$ is a strong support vertex of $T$. We note that $\pc(T) = r - 1 + s$. For each vertex of $V(T')$, select exactly one of its leaf-neighbors and let $S$ denote the resulting set of $s + 1$ leaves. The set $V(T) \setminus S$ is a TF-set of $T$, and so $F_t(T) \le |V(T)| - |S| = r + 2s$. Conversely by Lemma~\ref{leafLem}, $F_t(T) \ge r + 2s$. Consequently, $F_t(T) = r + 2s$. Therefore, letting $\ell = r - 1 + s$ and $k = s+1$, we note that $\ell \ge 3$ and $k \in [\ell - 1] \setminus \{1\}$, and that $\pc(T) = \ell$ and $F_t(T) = k + \ell$.~\smallqed

\medskip
We are now in a position to prove Theorem~\ref{t:main1}. Recall its statement.

\medskip
\noindent \textbf{Theorem~\ref{t:main1}}. \emph{If $T$ is a nontrivial tree, then $\pc(T) + 1 \le F_t(T) \le 2\pc(T)$, and all possible values of $F_t(T)$ in this range are possible. Further, the following hold. \\[-24pt]
\begin{enumerate}
\item $F_t(T) = \pc(T) + 1$ if and only if $\trim(T) = P_2$ or $\trim(T) \cong K_{1,n-1}$ for some $n \ge 4$.
\item $F_t(T) = 2\pc(T)$ if and only if $T$ has a unique minimum path cover and every path in this cover starts and ends at distinct leaves of~$T$.
\end{enumerate} }

\proof As an immediate consequence of Observation~\ref{o:FtZ},  Theorem~\ref{t:lowerbd}  and Theorem~\ref{t:cover}, if $T$ is a nontrivial tree, then $\pc(T) + 1 \le F_t(T) \le 2\pc(T)$. Further, $F_t(T) = \pc(T) + 1$ if and only if $\trim(T) = P_2$ or $\trim(T) \cong K_{1,n-1}$ for some $n \ge 4$. By Proposition~\ref{p:range}, all possible values of $F_t(T)$ in the range from $\pc(T) + 1$ to $2\pc(T)$ are possible. By Lemma~\ref{l:good_cover}, if $T$ has a unique minimum path cover and every path in this cover starts and ends at distinct leaves of~$T$, then $F_t(T) = 2\pc(T)$.

To complete the proof of Theorem~\ref{t:main1}, it suffices for us to prove that if $T$ is a nontrivial tree satisfying $F_t(T) = 2\pc(T)$, then $T$ has a unique minimum path cover and every path in this cover starts and ends at distinct leaves of~$T$. We proceed by induction on the order~$n \ge 2$ of a tree $T$ satisfying $F_t(T) = 2\pc(T)$. If $n = 2$, then $T \cong P_2$ and the result is immediate. This establishes the base case. Let $n \ge 3$ and assume that if $T'$ is a tree of order~$n'$ where $2 \le n' < n$ satisfying $F_t(T') = 2\pc(T')$, then $T'$ has a unique minimum path cover and every path in this cover starts and ends at distinct leaves of~$T'$. Let $T$ be a tree of order~$n$ satisfying $F_t(T) = 2\pc(T)$.

Suppose that $T \ne \trim(T)$. Let $T' = \trim(T)$. By supposition, $T'$ is a non-trivial tree of order less than~$n$. By Lemma~\ref{l:trimT} and Lemma~\ref{l:trimTb}, $F_t(T) = F_t(T')$ and $\pc(T) =\pc(T')$, implying that $F_t(T') = 2\pc(T')$. Applying the inductive hypothesis to $T'$, the tree $T'$ has a unique minimum path cover and every path in this cover starts and ends at distinct leaves of~$T'$. Thus, by Lemma~\ref{l:trimTb}, the tree $T$ has a unique minimum path cover and every path in this cover starts and ends at distinct leaves of~$T$.  Hence, we may assume that $T = \trim(T)$, for otherwise the desired result follows. With this assumption, we note that every edge in $T$ is incident with a vertex of degree at least~$3$. Thus, $T$ is not a path and every support vertex in $T$ has degree at least~$3$. In particular, $n \ge 4$.  We proceed further with the following series of claims.

\begin{claim}
\label{claim1}
Every strong support vertex in $T$ has exactly two leaf neighbors.
\end{claim}
\proof Suppose, to the contrary, that $T$ has a strong support vertex $v$ with three or more leaf neighbors. Let $v'$ be a leaf neighbor of $v$ in $T$ and consider the tree $T' = T - v'$. Let $S$ be a minimum TF-set of $T$. By Lemma~\ref{leafLem}, the set $S$ contains the vertex $v$ and all except possibly one leaf neighbor of~$v$. Renaming the leaf $v'$ if necessary, we may choose the set $S$ so that $v' \in S$. Thus since $S$ is a TF-set of $T$, the set $S \setminus \{v'\}$ is a TF-set of $T'$, and so $F_t(T') \le |S_v| - 1 = F_t(T) - 1$. Conversely, every minimum TF-set of $T'$ can be extended to a TF-set of $T$ by adding to it the vertex $v'$, and so $F_t(T) \le F_t(T') + 1$. Consequently, $F_t(T') = F_t(T) - 1$.

We show next $\pc(T') = \pc(T) - 1$. Every path cover in $T'$ can be extended to a path cover in $T$ by adding to it the trivial path consisting of the vertex $v'$, implying that $\pc(T) \le \pc(T') + 1$. To prove the reverse inequality, let $v_1$ and $v_2$ be two leaf neighbors of $v$ different from $v'$. By Lemma~\ref{l:strong}, there exists a minimum path cover, $\cP$ say, in $T$ that contains the path $v_1vv_2$. Let $P'$ be the path in $\cP$ that contains the vertex $v'$. Necessarily, the path $P'$ is the trivial path consisting of the vertex $v'$. Thus, $\cP \setminus \{P'\}$ is a path cover in $T'$, implying that $\pc(T') \le |\cP| - 1 = \pc(T) - 1$. Consequently, $\pc(T) = \pc(T') + 1$. Therefore by our earlier observations, $F_t(T) = F_t(T') + 1 \le 2\pc(T') + 1 = 2\pc(T) - 1$, contradicting our supposition that $F_t(T) = 2\pc(T)$. This completes the proof of Claim~\ref{claim1}.~\smallqed

\begin{claim}
\label{claim2}
If $\diam(T) \le 3$, then $T$ has a unique minimum path cover and every path in this cover starts and ends at distinct leaves of~$T$.
\end{claim}
\proof Suppose that $\diam(T) \le 3$. If $\diam(T) = 2$, then $T \cong K_{1,n-1}$ is a star. In this case, $\pc(T) = n-2$ and,
by Observation~\ref{simpleFormula}(d), $F_t(T) = n - 1$. Thus since $n \ge 4$,  $F_t(T) < 2\pc(T)$, a contradiction. Hence, $\diam(T) = 3$, implying that $T \cong S(r,s)$ is a double star. Since $T = \trim(T)$, we note that $r \ge 2$ and $s \ge 2$. Thus, $\pc(T) = s + t - 2$ and by Lemma~\ref{leafLem}, $F_t(T) = s + t$. Let $u$ and $v$ denote the two central vertices of the double star $T$. We note that $u$ and $v$ are the two (adjacent) vertices in $T$ that are not leaves. If $s + t \ge 5$, then $F_t(T) < 2\pc(T)$, a contradiction. Hence, $s + t = 4$, implying that $T \cong S(2,2)$ and that $T$ has a unique minimum path cover consisting of two paths, namely a path containing $u$ and its two leaf neighbors and a path containing $v$ and its two leaf neighbors.~\smallqed

\medskip
By Claim~\ref{claim2}, we may assume that $\diam(T) \ge 4$, for otherwise the desired result follows. Let $u$ and $r$ be two vertices at maximum distance apart in $T$. Necessarily, $u$ and $r$ are leaves and $d(u,r) = \diam(T)$. We now root the tree $T$ at the vertex $r$. Let $v$ be the parent of $u$, $w$ the parent of $v$, $x$ be the parent of $w$, and $y$ the parent of~$x$.
We note that if $\diam(T) = 4$, then $y = r$; otherwise, $y \ne r$. By our earlier assumptions, every support vertex has degree at least~$3$. In particular, $d_T(v) \ge 3$, and so $v$ is a strong support vertex. By Claim~\ref{claim1}, $d_T(v) = 3$. Let $u_1$ and $u_2$ be the two leaf neighbors of $v$, where $u = u_1$.

Let $T'$ be the tree obtained from $T$ by deleting $v$ and its two children; that is, $T' = T - \{v,u_1,u_2\}$. Let $T'$ have order~$n'$, and so $n' = n - 3$. Since $\diam(T) \ge 4$, we note that $n' \ge 3$. Let $P'$ be the path $u_1vu_2$. By Lemma~\ref{l:strong}, there exists a minimum path cover, $\cP$ say, in $T$ that contains the path $P'$. Since $\cP \setminus \{P'\}$ is a path cover in $T'$, we note that $\pc(T') \le |\cP| - 1 = \pc(T) - 1$. Every minimum TF-set of $T'$ can be extended to a TF-set of $T$ by adding to it the vertices $v$ and $u_1$, implying that $F_t(T) \le F_t(T') + 2$. Therefore by our earlier observations,
\begin{equation}
\label{Eq1}
2\pc(T) = F_t(T) \le F_t(T') + 2 \le 2\pc(T') + 2 \le 2(\pc(T) - 1) + 2 = 2\pc(T).
\end{equation}

Hence we must have equality throughout the above Inequality Chain~(\ref{Eq1}), implying that $F_t(T') = 2\pc(T')$ and $\pc(T') = \pc(T) - 1$. Applying the inductive hypothesis to $T'$, the tree $T'$ has a unique minimum path cover $\cP'$ and every path in $\cP'$ starts and ends at distinct leaves of~$T'$. Let $\cP' = \{Q_1,\ldots,Q_k\}$, where $Q_1$ is the path that contains the vertex~$w$.

\begin{claim}
\label{claim3}
$d_T(w) \ge 3$.
\end{claim}
\proof Suppose, to the contrary, that $d_T(w) = 2$, implying that $w$ is a leaf in $T'$ with the vertex~$x$ as its neighbor. Let $T^* = T - \{u_1,u_2\}$. Let $Q_1^*$ be the path obtained from $Q_1$ by adding to it the vertex $v$ and the edge $vw$. If $k \ge 2$, let  $Q_i^* = Q_i$ for $i \in [k] \setminus \{1\}$. Let $\cQ^* = \{Q_1^*,\ldots,Q_k^*\}$. Since $\cP'$ is a unique minimum path cover in $T'$, we note that $\cP^*$ is a unique minimum path cover in $T^*$. In particular, $\pc(T') = \pc(T^*)$. Further since every path in $\cP'$ starts and ends at distinct leaves of~$T'$, every path in $\cP^*$ starts and ends at distinct leaves of~$T^*$. Let $S^*$ consist of a leaf and it neighbor from every path in $\cQ^*$, where we choose $S^*$ so that $\{v,w\} \subseteq S^*$. By Lemma~\ref{l:good_cover}, the set $S^*$ is a TF-set of $T^*$. The set $S^*$ can be extended to a TF-set of $T$ by adding to it the vertex $u$, implying that $F_t(T) \le |S^*| + 1 = 2\pc(T^*) + 1 = 2\pc(T') + 1 = 2(\pc(T) - 1) + 1 < 2\pc(T)$, a contradiction. Therefore, $d_T(w) \ge 3$.~\smallqed

\medskip
By Claim~\ref{claim3}, $d_T(w) \ge 3$, implying that $w$ is not a leaf in $T'$. The vertex $w$ is therefore an internal vertex on the path $Q_1$. Recall that $P'$ is the path $u_1vu_2$. We now consider that path cover $\cP = \cP' \cup \{P'\} = \{P',Q_1,\ldots,Q_k\}$. As observed earlier, $\pc(T') = \pc(T) - 1$, and so $|\cP| = |\cP'| + 1 = \pc(T') + 1 = \pc(T)$, and so $\cP$ is a minimum path cover in $T$.

\begin{claim}
\label{claim4}
$\cP$ is the unique minimum path cover in $T$.
\end{claim}
\proof Suppose, to the contrary, that there is a minimum path cover, $\cP^*$ say, that is different from $\cP$. If $P'$ is a path in $\cP^*$, then $\cP^* \setminus \{P'\}$ is a minimum path cover in $T'$ different from $\cP'$, a contradiction. Hence, $P'$ is not a path in $\cP^*$. Let $P_v^*$ be the path in $\cP^*$ that contains the vertex $v$, and so $P_v^* \ne P'$. By the minimality of the path cover $\cP^*$, exactly one of $u_1$ and $u_2$, say $u_2$, belong to the path $P_v^*$. Thus, the vertex $u_1$ belong to a trivial path, say $P_u^*$, in $\cP^*$ consisting only of the vertex $u_1$. If $P_v^*$ does not contain the vertex $w$, then $\cP^* \setminus \{P_u^*,P_v^*\}$ is a path cover in $T'$ of size~$|\cP^*| - 2 = \pc(T) - 2 = \pc(T') - 1$, a contradiction. Therefore, $P_v^*$ contains the vertex $w$. Let $P_w^*$ be obtained from $P_v^*$ by deleting from it the vertices $v$ and $u_2$, and so $P_w^* = P_v^* - \{v,u_2\}$. We now consider the path cover of $T'$ consisting of the path $P_w^*$ together with all paths in $\cP^*$ different from  $P_u^*$ and $P_v^*$. The resulting path cover in $T'$ has size $|\cP^*| - 1 = \pc(T) - 1 = \pc(T')$ and is therefore a minimum path cover in $T'$. However, the path $P_w^*$ in this path cover has as one of its end the vertex $w$, which is not a leaf in $T'$, and this path cover is therefore different from $\cP'$. This contradicts the fact that $\cP'$ is the unique minimum path cover in $T'$. Hence, $\cP$ is the unique minimum path cover in $T$.~\smallqed

\medskip
By Claim~\ref{claim4}, the tree $T$ has a unique minimum path cover, namely $\cP$. By our earlier observations, every path in $\cP$ starts and ends at distinct leaves of~$T$. This completes the proof of Theorem~\ref{t:main1}.~\qed

\section{Proof of Theorem~\ref{t:main2}}
\label{S:main2}

In this section we prove Theorem~\ref{t:main2}. First we present the following lemma showing that every tree $T$ in the family~$\cT$ satisfies $F_t(T) = \alpha'(T) + \pc(T)$.

\begin{lem}
\label{l:match}
If $T \in \cT$, then $F_t(T) = \alpha'(T) + \pc(T)$.
\end{lem}
\proof We proceed by induction on the order~$n \ge 2$ of a tree $T \in \cT$. If $n \in \{2,3\}$, then $T = P_2$ or $T = P_3$. In both cases, the result is immediate noting that $F_t(T) = 2$ and $\alpha'(T) = \pc(T) = 1$. This establishes the base cases. Let $n \ge 4$ and assume that if $T' \in \cT$ is a tree of order~$n'$ where $n' < n$, then $F_t(T') = \alpha'(T') + \pc(T')$. Let $T \in \cT$ be a tree of order~$n$.

Let $F$ be the underlying tree of $T \in \cT$ and let $A$ be the attacher set of $T$. Let $B = V(F) \setminus A$. Thus, $A \subseteq V(F)$ and either $A = V(F)$ or $A \subset V(T)$ and $B$ is an independent set in $F$ containing no leaf of $F$. Further, the set of support vertices in $T$ is precisely the set of attacher vertices (that belong to $A$), and each attacher vertex is a strong support vertex of $T$ with all its leaf neighbors incident with pendant edges that were added to $F$ when forming $T$.

Suppose that $T$ contains a support vertex $v$ with three or more leaf neighbors. Let $v'$ be a leaf neighbor of $v$ in $T$ and consider the tree $T' = T - v'$. We note that $T' \in \cT$ and that $T$ and $T'$ have the identical underlying tree, namely $F$, and the same attacher set $A$. Applying the inductive hypothesis to $T'$, the tree $T'$ satisfies $F_t(T') = \alpha'(T') + \pc(T')$. We note that $\alpha'(T') = \alpha'(T)$. Identical arguments as in the proof of Claim~\ref{claim1} of Theorem~\ref{t:main1} show that $F_t(T') = F_t(T) - 1$ and $\pc(T') = \pc(T) - 1$. Thus, $F_t(T) = F_t(T') + 1 = \alpha'(T') + \pc(T') + 1 = \alpha'(T) + \pc(T)$. Hence, we may assume that every support vertex in $T$ has exactly two leaf neighbors, for otherwise the desired result follows.

We show that $\alpha(T) = |A|$. Let $M$ be a maximum matching in $T$. By the maximality of $M$, each attacher vertex of $T$ is incident with an edge of $M$. Let $v \in A$ denote an arbitrary attacher vertex of $T$ and let $v'$ denote one of its leaf neighbors. If $vv' \notin M$, then we can simply replace the edge of $M$ incident with $v$ with the edge $vv'$. Hence, we may assume that $vv' \in M$. More generally, we can choose $M$ to contain $|A|$ pendant edges in $T$ associated with the $|A|$ attacher vertices in $A$. With this choice of $M$, we note that a leaf that is not incident with one of these $|A|$ pendant edges does not belong to $M$. Thus, the only possibly additional edges in $M$ are edges with both ends in $F$. If $A \subset V(T)$, then noting that $B$ is an independent set in $F$ and the only neighbors in $T$ of vertices in $B$ are attacher vertices in $A$ which are already matched under $M$ with one of their leaf neighbors, no vertices of $B$ are incident with an edge of $M$. This implies that $\alpha(T) = |M| = |A|$.

We show next that $\pc(T) = |V(F)|$. For each attacher vertex $v \in A$ in $T$, let $v_1$ and $v_2$ denote its two leaf neighbors. By our earlier assumption, all other neighbors of $v$ belong to the underlying tree $F$. Let $\cP$ be a minimum path cover in $T$. By an identical proof as shown in Lemma~\ref{l:strong} we can choose $\cP$ so that it contains the path $v_1vv_2$ for every such attacher vertex $v$. As observed earlier, if $A \subset V(T)$, then $B$ is an independent set in $F$ and the only neighbors in $T$ of vertices in $B$ are attacher vertices in $A$. Thus, each vertex in $B$ belongs to a path in $\cP$ that is a trivial path consisting only of that vertex. Thus, each vertex in $V(F)$ belongs to a distinct path in the path cover $\cP$, implying that $\pc(T) = |\cP| = |V(F)|$.

Finally, we show that $F_t(T) = |V(F)| + |A|$. Among all minimum TF-set of $T$, let $S$ be chosen to contain as few leaves as possible. As observed earlier, each attacher vertex, $v$ say, of $T$ is a strong support vertex with two leaf neighbors, say $v_1$ and $v_2$. By Lemma~\ref{leafLem}, the set $S$ contains the vertex $v$ and at least one of $v_1$ and $v_2$. If both $v_1$ and $v_2$ belong to $S$, then by the minimality of the TF-set $S$, there is a neighbor $v'$ of $v$ not in $S$. Such a vertex necessarily belongs to the set $B$. Replacing the vertex $v_1$ in $S$ with the vertex $v'$ produces a new minimum TF-set of $T$ that contains fewer leaves than does the set $S$, a contradiction. Hence exactly one leaf neighbor of every attacher vertex does not belong to $S$. We show next that every vertex in $B$ belongs to $S$. If this is not the case, then let $w$ be a vertex in $B$ that does not belong to $S$. As observed earlier, every neighbor of $w$ in $T$ is an attacher vertex (that belongs to $A$) with one of its leaf neighbors not in $S$. This implies, however, that $S$ is not a forcing set since the vertex $w$ cannot be colored in the forcing process starting with the set $S$, a contradiction. Hence, $B \subset S$, implying that $V(F) \subset S$ and that exactly one leaf neighbor of every vertex in $A$ belongs to $S$. Thus, $F_t(T) = |S| = |V(F)| + |A|$. As observer earlier, $\alpha(T) = |A|$ and $\pc(T) = |V(F)|$. Therefore, $F_t(T) = \alpha'(T) + \pc(T)$.~\qed

\medskip 
We are now in a position to prove Theorem~\ref{t:main2}. Recall its statement.

\medskip
\noindent \textbf{Theorem~\ref{t:main2}}. \emph{If $T$ is a nontrivial tree, then  $F_t(T) \le \alpha'(T) + \pc(T)$, with equality if and only if $T \in \cT$.}

\proof  We proceed by induction on the order~$n \ge 2$ of a tree $T$ to show that $F_t(T) \le \alpha'(T) + \pc(T)$ and that if equality holds, then $T \in \cT$. If $n \in \{2,3\}$, then $T = P_2$ or $T = P_3$. In both cases, the result is immediate noting that $F_t(T) = 2$ and $\alpha'(T) = \pc(T) = 1$, and $T \in \cT$. This establishes the base cases. Let $n \ge 4$ and assume that if $T'$ is a tree of order~$n'$ where $2 \le n' < n$, then $F_t(T') \le \alpha'(T') + \pc(T')$, with equality if and only if $T' \in \cT$. Let $T$ be a tree of order~$n$. If $T \cong P_n$ is a path, then $F_t(T) = 2$ and $\pc(T) = 1$. However since $n \ge 4$, we note that in this case $\alpha'(T) \ge 2$, and so $F_t(T) < \alpha'(T) + \pc(T)$. Hence, we may assume that $T$ is not a path, for otherwise the desired result follows.

\begin{unnumbered}{Claim~A}
If $T$ has a support vertex with three of more leaf neighbors, then the desired result follows.
\end{unnumbered}
\proof  Suppose that $T$ has a support vertex $v$ with three or more leaf neighbors. Let $v'$ be a leaf neighbor of $v$ in $T$ and consider the tree $T' = T - v'$. We note that $\alpha'(T') = \alpha'(T)$. Identical arguments as in the proof of Claim~\ref{claim1} of Theorem~\ref{t:main1} show that $F_t(T') = F_t(T) - 1$ and $\pc(T') = \pc(T) - 1$.  Applying the inductive hypothesis to $T'$, we therefore have that
\begin{equation}
\label{Eq2}
F_t(T) = F_t(T') + 1 \le (\alpha'(T') + \pc(T')) + 1 = \alpha'(T) + \pc(T).
\end{equation}

Further, suppose that $F_t(T) = \alpha'(T) + \pc(T)$. In this case, we must have equality throughout the above Inequality Chain~(\ref{Eq2}). Thus, $F_t(T') = \alpha'(T') + \pc(T')$, and so by the inductive hypothesis, $T' \in \cT$. We note that the vertex $v$ is a strong support vertex of $T'$, implying that the vertex $v$ is a vertex of the underlying tree used to construct $T' \in \cT$ and the leaf neighbors of $v$ do not belong to the underlying tree. This in turn implies that $T \in \cT$ (and that both $T'$ and $T$ have the same underlying tree). This completes the proof of Claim~A.~\smallqed

\medskip
By Claim~A, we may assume that every strong support vertex in $T$ has exactly two leaf neighbors, for otherwise the desired result follows.

\begin{unnumbered}{Claim~B}
If $T \ne \trim(T)$, then $F_t(T) < \alpha'(T) + \pc(T)$.
\end{unnumbered}
\proof Suppose that $T \ne \trim(T)$. Let $T' = \trim(T)$. By supposition, $T'$ is a non-trivial tree of order less than~$n$. By Lemma~\ref{l:trimT} and Lemma~\ref{l:trimTb}, $F_t(T) = F_t(T')$ and $\pc(T) = \pc(T')$. Contracting edges cannot increase the matching number, implying that $\alpha'(T') \le \alpha'(T)$.  Applying the inductive hypothesis to $T'$, we therefore have that
\begin{equation}
\label{Eq3}
F_t(T) = F_t(T') \le \alpha'(T') + \pc(T') \le \alpha'(T) + \pc(T).
\end{equation}

We show next that $F_t(T) < \alpha'(T) + \pc(T)$. Suppose to the contrary that $F_t(T) = \alpha'(T) + \pc(T)$. In this case, we must have equality throughout the above Inequality Chain~(\ref{Eq3}). Thus, $\alpha'(T') = \alpha(T)$. Further, $F_t(T') = \alpha'(T') + \pc(T')$, and so by the inductive hypothesis, the tree $T' \in \cT$. Let $F$ be the underlying tree of $T' \in \cT$ and let $A'$ be the attacher set of $T'$. Let $B' = V(F) \setminus A'$. Thus, $A' \subseteq V(F)$ and either $A' = V(F)$ or $A' \subset V(F)$ and $B'$ is an independent set in $F$ containing no leaf of $F$. Further, the set of support vertices in $T'$ is precisely the set of attacher vertices (that belong to $A'$), and each attacher vertex is a strong support vertex of $T'$ with all its leaf neighbors incident with pendant edges that were added to $F$ when forming $T'$. By construction of trees in the family~$\cT$, every support vertex of $T' \in \cT$ is a strong support vertex.

By definition of a trimmed tree, the tree $T$ can be rebuilt from the tree $T' = \trim(T)$ by subdividing edges of $T'$. Subdividing edges cannot decrease the matching number. However as observed earlier, $\alpha'(T') = \alpha(T)$, implying that at every stage of the rebuilding process starting from $T'$, whenever we subdivide an edge the matching number remains unchanged. We show, however, that this is not the case. Let $e$ be the first edge of $T' = \trim(T)$ that is subdivided in this reconstruction process to rebuild the tree $T$, and let $T^*$ be obtained from $T'$ by subdividing the edge $e$.

Suppose firstly that $e$ is a pendant edge of $T'$, say $e = vv_1$ where $v_1$ is a leaf of $T'$. Thus, $v$ is the strong support vertex of $T'$, or, equivalently, $v$ is an attacher vertex of $T'$, and so $v \in A'$. Let $u$ be the new vertex of degree~$2$ resulting from subdividing the edge~$e$, and so $u$ is a support vertex of $T^*$ with $v_1$ as its leaf neighbor and $v$ as its non-leaf neighbor. Let $v_2$ be a leaf neighbor of $v$ in $T'$ different from $v_1$. Let $M$ be a maximum matching in $T'$. By the maximality of $M$, the vertex $v$ is incident with an edge of $M$. If $vv_2 \notin M$, then we can simply replace the edge of $M$ incident with $v$ with the edge $vv_2$. Hence, we may assume that $vv_2 \in M$. But then $M \cup \{uv_1\}$ is a matching in $T^*$, implying that $\alpha'(T^*) > |M| = \alpha'(T')$. Since $\alpha'(T) \ge \alpha'(T^*)$, this implies that  $\alpha(T) > \alpha'(T')$, a contradiction.

Suppose next that $e = uv$ is not a pendant edge of $T'$. Thus, the edge $e$ belongs to the underlying tree $F$ of $T'$. By definition of a trimmed tree, the edge $e$ is incident with a vertex of degree~$2$ and with a vertex of degree at least~$3$. Renaming $u$ and $v$ if necessary, we may assume that $u$ has degree~$2$ and $v$ has degree at least~$3$ in $T'$. Let $w$ be the neighbor of $u$ different from $v$. Since each attacher vertex of $T'$ has degree at least~$3$ and since the set $B'$ is an independent set in $T'$, this implies that $u \in B'$ and $\{v,w\} \subseteq A'$. Thus, both $v$ and $w$ are strong support vertices in $T'$ with their leaf neighbors outside $F$.

Let $u'$ be the new vertex of degree~$2$ resulting from subdividing the edge~$e = uv$, and so $u'$ has as its neighbors in $T^*$ the vertices $u$ and $v$. Let $M$ be a maximum matching in $T'$. By the maximality of $M$, both vertices $v$ and $w$ are incident with edges of $M$. Let $v'$ and $w'$ be arbitrary leaf neighbors of $v$ and $w$, respectively, in $T'$. If $vv' \notin M$, then we can simply replace the edge of $M$ incident with $v$ with the edge $vv'$. Hence, we may assume that $vv' \in M$. Analogously, we may assume that $ww' \in M$. But then $M \cup \{uu'\}$ is a matching in $T^*$, implying that $\alpha'(T^*) > |M| = \alpha'(T')$. Since $\alpha'(T) \ge \alpha'(T^*)$, this implies that $\alpha(T) > \alpha'(T')$, a contradiction. Therefore, $F_t(T) < \alpha'(T) + \pc(T)$.  This completes the proof of Claim~B.~\smallqed

\medskip
By Claim~B, we may assume that $T = \trim(T)$, for otherwise the desired result follows. With this assumption, we note that every edge in $T$ is incident with a vertex of degree at least~$3$. In particular, every support vertex in $T$ has degree at least~$3$. By our earlier assumption, every strong support vertex in $T$ has exactly two leaf neighbors. Thus since $n \ge 4$, we note that $T$ is not a star. Hence, $\diam(T) \ge 3$. Let $u$ and $r$ be two vertices at maximum distance apart in $T$. Necessarily, $u$ and $r$ are leaves and $d(u,r) = \diam(T)$. We now root the tree $T$ at the vertex $r$. Let $v$ be the parent of $u$, $w$ the parent of $v$, and $x$ be the parent of $w$. Our earlier assumptions imply that $d_T(v) =3$. Let $u_1$ and $u_2$ be the two children of $v$, where $u = u_1$.

Let $T'$ be the tree obtained from $T$ by deleting $v$ and its two children; that is, $T' = T - \{u_1,u_2,v\}$. Let $T'$ have order~$n'$, and so $n' = n - 3$. Since $\diam(T) \ge 3$ and $T = \trim(T)$, we note that $n' \ge 3$. Let $P'$ be the path $u_1vu_2$. By Lemma~\ref{l:strong}, there exists a minimum path cover, $\cP$ say, in $T$ that contains the path $P'$. Thus, $\cP' = \cP \setminus \{P'\}$ is a path cover in $T'$, implying that $\pc(T') \le |\cP| - 1 = \pc(T)  - 1$. Further, we note that $\alpha'(T') = \alpha'(T) - 1$. Every minimum TF-set of $T'$ can be extended to a TF-set of $T$ by adding to it the vertices $u_1$ and $v$, implying that $F_t(T) \le F_t(T') + 2$. Therefore by our earlier observations,
\[
\begin{array}{lcl}
F_t(T) & \le & F_t(T') + 2 \1 \\
& \le & \alpha'(T') + \pc(T') + 2 \1 \\
& \le & (\alpha(T) - 1) + (\pc(T)  - 1) + 2 \1 \\
& = & \alpha(T) + \pc(T).
\end{array}
\]

This establishes the desired upper bound. Suppose next that $F_t(T) = \alpha'(T) + \pc(T)$. In this case, we must have equality throughout the above inequality chain. Thus, $\pc(T') = \pc(T) - 1$, $F_t(T) = F_t(T') + 2$ and $F_t(T') = \alpha'(T') + \pc(T')$. Applying the inductive hypothesis to $T'$, the tree  $T' \in \cT$. Let $F$ be the underlying tree of $T' \in \cT$ and let $A'$ be the attacher set of $T'$. We now consider two possibilities, depending on whether $w \in V(F)$ or $w \notin V(F)$.

Suppose firstly that $w$ does not belong to the underlying tree $F$ of $T'$, implying that $w$ is a leaf of $T'$ and that $x$ be the (unique) neighbor of $w$ in $T'$. We note that $x \in A'$ and that by our earlier assumptions, the vertex $x$ is either a support vertex in $T'$ with exactly two leaf neighbors or a support vertex in $T'$ with exactly three leaf neighbors.

We show that the vertex $x$ has exactly three leaf neighbors in $T'$. Suppose, to the contrary, that $x$ has exactly two leaf neighbors in $T'$. Let $L'$ be a set of $|A'|$ leaves in $T'$ consisting of exactly one leaf neighbor of every vertex of $A'$ in $T'$. Further, we choose $L'$ so that $w \in L$. We note that the set $V(F) \cup L$ is a minimum TF-set in $T'$, and so $F_t(T') = |V(F)| + |L| = |V(F)| + |A'| = n' - |A'|$. However, the set $(V(F) \setminus \{z\}) \cup L \cup \{v,u_1\}$ is a TF-set of $T$, where as the first vertex played in the forcing process we play the vertex $w$ (of degree~$2$ in $T$) which forces the vertex $z$ to be colored, as the second vertex we play the vertex $v$ which colors the vertex $u_2$, and thereafter we play the identical sequence of vertices in the forcing process in $T'$ starting with the set $V(F) \cup L$ that results in all $V(T')$ colored. Thus, $F_t(T) \le |V(F)| - 1 + |L| + 2 = |V(F)| + |A'| + 1 = n' - |A'| + 1 = F_t(T') + 1$, contradicting our earlier observation that $F_t(T) = F_t(T') + 2$.

Hence, the vertex $z$ has exactly three leaf neighbors in $T'$. This implies that $T \in \cT$, where we note that the underlying tree $U$ of $T$ is obtained from the tree $F$ by adding to it the vertices $v$ and $w$ and the edges $vw$ and $wz$, and where the attacher set $A$ of $T$ is the set $A = A' \cup \{v\}$.

Suppose secondly that the vertex $w$ belongs to the underlying tree $F$. This implies that $T \in \cT$, where we note that the underlying tree $U$ of $T$ is obtained from the tree $F$ by adding to it the vertex $v$ and the edge $vw$, and where the attacher set $A$ of $T$ is the set $A = A' \cup \{v\}$. This completes the proof that if $T$ is a nontrivial tree, then  $F_t(T) \le \alpha'(T) + \pc(T)$ and if equality holds, then $T \in \cT$. By Lemma~\ref{l:match}, if $T \in \cT$, then $F_t(T) = \alpha'(T) + \pc(T)$. This completes the proof of Theorem~\ref{t:main2}.~\qed

\section{Closing Remarks}
\label{S:remark}

The total forcing number of a tree $T$ and its path cover number are related by the inequality chain $\pc(T) + 1 \le F_t(T) \le 2\pc(T)$. In this paper, we characterize the extremal trees achieving equality in these bounds (see, Theorem~\ref{t:main1}). We remark that the inequality $F_t(G) \le 2\pc(G)$ is not true for general graphs $G$. Even for the class of cubic graphs, there is no constant $C$ such that $F_t(G) \le C \times \pc(G)$ holds for every connected cubic graph $G$.

Our second main result, namely Theorem~\ref{t:main2}, shows that the total forcing number of a tree $T$ is related to its matching number and path cover number by the inequality $F_t(T) \le \alpha'(T) + \pc(T)$. Further, we characterize the trees achieving equality in this bound. We remark that the inequality $F_t(G) \le \alpha'(G) + \pc(G)$ is not true for general graphs $G$. As simple counterexamples, take $G = K_n$ where $n \ge 5$ or $G = K_{k,k}$ where $k \ge 4$.

\end{document}